\documentclass[12pt]{article} \textwidth=6in \oddsidemargin=0in

\usepackage{pst-all}
\usepackage[dvips]{graphics}

\usepackage{amssymb}
\usepackage{color} 
\begin{document}
\def\ov{\over} \def\be{\begin{equation}}
\def\ee{\end{equation}} \def\x{\xi} \def\s{\sigma} \def\iy{\infty}
\def\bc{\begin{center}} \def\ec{\end{center}}
\def\Z{\mathbb Z} \def\inv{^{-1}}  \def\C{\mathcal C} 
\def\ep{\varepsilon} \def\X{\mathcal X} \def\g{\gamma} 
\def\d{\delta} \def\Phat{\widehat P} \def\d{\delta} 
\def\R{\mathbb R} \def\ph{\varphi} \def\sq{\sqrt{2 s}} 
\def\g{\gamma} \newcommand{\xs}[1]{\x_{\s(#1)}} \def\B{\mathbb B}
\newcommand{\vs}[1]{v_{\s(#1)}} \def\sqq{{\sqrt{2\ov s}}} 
\def\ps{\psi} \def\cd{\cdots} \def\Bn{\B_{n-1}}
\def\ld{\ldots} \def\S{\mathbb S} \def\Wn{\mathbb W^{n-1}} 
\def\1{{\bf 1}} \def\Rn{\R^{n-1}} \def\Sn{\S_{n-1}} \def\e{\eta}
\def\a{\alpha} \def\P{{\mathcal P}} \def\l{\ell} \def\t{\tau}
\def\lb{\langle} \def\rb{\rangle} \def\pst{\widetilde\ps} 
\def\sech{{\rm sech}\,}
\def\1{{\bf 1}} \def\a{\alpha} \def\b{\beta} \def\m{\mu} 
\def\Ph{\Phi} \def\Phh{\widehat\Ph} \def\r{\rho} \def\noi{\noindent}
\def\sp{\vskip2ex} \def\bx{{\bf x}} \def\by{{\bf y}} \def\bz{{\bf z}} \def\Q{\mathcal Q} \def\p{\mathfrak p} \def\ph{\widehat\p}
\def\Ps{\Psi} \def\Psh{\widehat\Ps} \def\Q{\mathcal Q} 
\def\Phat{\widehat \gP}\def\L{\mathcal L} 
\def\LP{\widehat\P} \def\Phat{\widehat P} \def\wh{\widehat}
\def\Kt{\tilde{K}} \def\Kh{\wh K} \def\D{\Delta} \def\E{\mathcal E}

\hfill April 1, 2013
\bc{\bf \large The Asymmetric Simple Exclusion Process\\ \vskip.5ex
 with an Open Boundary}\ec

\bc{\large\bf Craig A.~Tracy}\\
{\it Department of Mathematics \\
University of California\\
Davis, CA 95616, USA}\ec
 
\bc{\large \bf Harold Widom}\\
{\it Department of Mathematics\\
University of California\\
Santa Cruz, CA 95064, USA}\ec

\bc{\bf I. Introduction}\ec

In previous work \cite{tw3} the authors considered the asymmetric simple exclusion process (ASEP) where particles are confined to the nonnegative integers $\Z^+=\{0,1,2,\ldots\}$. Each particle waits exponential time, and then with probability $p$ it moves one step to the right if the site is unoccupied, otherwise it does not move; and with probability $q=1-p$ a particle not at 0 moves one step to the left if the site is unoccupied, otherwise it does not move. For $n$-particle ASEP a possible cofigurations is 
\[\bx=\{x_1,\ld,x_n\},\quad (0\le x_1<\cd<x_n).\]
The $x_i$ are the occupied sites. We denote by $\X_n$ the set of possible configurations for $n$-particle ASEP, and by $\p_n(\bx,\by;t)$ the probability that at time~$t$ the system is in configuration~$\bx$ given that initially it was in configuration $\by$. (We shall drop the subscript ``$n$'' when it is understood.)

In \cite{tw3} a formula was found for this probability. It was the sum over the Weyl group $\B_n$ of multiple integrals. (For ASEP on $\Z$ it was a sum over the permutation group $\S_n$ \cite{tw1,tw2}.)

Here we consider the ASEP on $\Z^+$ with an open boundary at zero. The stationary  measure for ASEP on the finite lattice $[1,L]$ or on the semi-infinite lattice $\Z^{+}$ with boundaries  connected
to reservoirs  has been the subject of much research starting with Derrida \textit{et al.}\ \cite{DEHP}.  We refer the reader to the recent work of Sasamoto and Williams \cite{SW} for   an up-to-date account of these developments.  Here we consider the time-dependent properties of ASEP on $\Z^+$ with an open boundary.  Specifically,
the point $0$ is connected to a reservoir where a particle 
is \textit{injected} into site $0$ from the reservoir at a rate $\a$, assuming that the site $0$ is empty, and a particle at site $0$ is \textit{ejected} into the reservoir at a rate $\b$. Now the number of particles is not conserved and for ASEP with open boundary the configuration $\bx$ may lie in $\X_n$ while $\by$ may lie in $\X_m$ with $m\ne n$.  

We find an infinite tri-diagonal matrix with operator entries in which the Laplace transforms of the probabilities can be read off from the entries of the inverse matrix. When either $\a=0$ or $\b=0$ the matrix is triangular and so the inverse can be computed more explicitly. The result is obtained by solving a system of differential equations for the probabilities. The final formulas involve inverses of operators with kernels the Laplace transforms of certain $\p(\bx,\by;t)$ obtained in \cite{tw3}. 

There are two special cases in which the results are more explicit. For TASEP with $p=1$, the inverse operator is computable in terms of $\p(\bx,\by;t)$ itself, and the probabilities are given in terms of certain determinants. For SSEP ($p=q$) and general $\a$ and $\b$ we find formulas analogous to the ones described above for the probability that sites $x_1,\ld,x_n$ are occupied. This is for infinite systems as well as finite ones. 

We state the formulas for $\p(\bx,\by;t)$ in the Appendix.

\bc{\bf II. Statement of results}\ec

We denote by $\E_n$ the Lebesgue space $L^1(\X_n)$. From the fact that  
\be \sum_{\bx\in\X_n}\p(\bx,\by;t)=1\label{psum}\ee
for each $\by$, it follows that the operator on $\E_n$ with this kernel is bounded with norm one. We denote the Laplace transform of $\p(\bx,\by;t)$ by $\ph(\bx,\by;s)$:
\[\ph(\bx,\by;s)=\int_0^\iy \p(\bx,\by;t)\,e^{-st}\,dt.\]
The operator with this kernel is bounded on $\E_n$ with norm at most $({\rm Re}\,s)\inv$. We denote it by $L_n(s)$.\footnote{When $n=0$ we interpret $\p(\bx,\by;t)$ as 1, and so $L_0(s)$ is multiplication by $s\inv$.} In the results stated below it is tacitly assumed that Re\,$s$ is sufficiently large. 

Now for ASEP with open boundary at zero, we define $P_n(\bx;\,t)$ to be the probability that the system is in configuration $\bx\in\X_n$ at time~$t$. (We shall usually drop the ``$\bx$'' in the notation, and do not specify an initial configuration.) We denote its Laplace transform by $\Phat_n(s)$. 

We define vector functions
\[\Phat(s)=(\Phat_n(s))_{n\ge0},\ \ \ \ P(0)=(P_n(0))_{n\ge0},\]
belonging to the direct sum $\sum_{n=0}^\iy \E_n$.
 
We define operators $A_n:\E_{n-1}\to \E_n$ and 
$B_n:\E_{n+1}\to \E_n$ by
\be(A_nF)(x_1,\ld,x_n)=\d(x_1)\,F(x_2,\ld,x_n),\label{A}\ee
\be(B_nF)(x_1,\ld,x_n)=(1-\d(x_1))\,F(0,x_1,\ld,x_n),\footnote{When $n=0$ we interpret $\d(x_1)$ as zero. In particular $A_0=0$ and $B_0\,F=F(0)$.}\label{B}\ee

Then we define matrices $\d,\,L(s),\,A,\, B$, with operator entries, acting on $\sum_{n=0}^\iy \E_n$. The first is diagonal with $(n,n)$-entry multiplication by $\d(x_1)$, the second is diagonal with $n,n$-entry $L_n(s)$, the third is subdiagonal (one diagonal below the main diagonal) with $n,n-1$-entry $A_n$, and the last is superdiagonal (one diagonal above the main diagonal) with $n,n+1$-entry $B_n$. 
\sp   

\noi{\bf Theorem 1}. With this notation we have
\be\Phat(s-\a)=\Big(I-L(s)\,((\a-\b)\,\d+\a\,A+\b\,B)\Big)\inv\,L(s)\,P(0).\label{Th1}\ee

There are expressions for the entries of the inverse operator as infinite series of products. But when either $\b=0$ or $\a=0$ the operator has only one subdiagonal or one superdiagonal and each entry of the inverse is a single product. We state the results as recursion formulas. Define 
\be M_n(s)=(I-(\a-\b)\,L_n(s)\,\d)\inv.\label{Mn}\ee

\noi{\bf Corollary 1.1}. Suppose $\b=0$ and that initially there are $k$ particles at $\by\in\X_k$. Then 
\[\Phat_k(s-\a)=M_{k}(s)\,L_{k}(s)\,\d_{\by},\footnote{This is interpreted as $s\inv$ when $k=0$.}\]
and when $n>k$
\[\Phat_n(s-\a)=\a\,M_{n}(s)\,L_{n}(s)\,A_n\,\Phat_{n-1}(s-\a).\]
\sp
 
\noi{\bf Corollary 1.2}. Suppose $\a=0$ and that initially there are $k$ particles at $\by\in\X_k$. Then 
\[\Phat_k(s)=M_{k}(s)\,L_{k}(s)\,\d_{\by},\]
and when $n<k$
\[\Phat_n(s)=\b\,M_{n}(s)\,L_{n}(s)\,B_n\,\Phat_{n+1}(s).\]
\sp

In connection with the corollaries we show the following.
\sp

\noi{\bf Remark 1.1}. The operators appearing in the inverses can be replaced by lower-dimensional ones. This will be useful for computation. Define 
\[\X_n^+=\{\{x_1,\ld,x_n\}\in\X_n:x_1>0\},\ \ \ \E_n^+=L^1(\X_n^+),\]
and then operators: 

\noi$L_{n-1}^0(s):\E_{n-1}^+\to\E_{n-1}^+$ with kernel $\ph((0,\bx),(0,\by);s)$, 

\noi  $L_{n,n-1}^0(s):\E_{n-1}^+\to\E_n$ with kernel $\ph(\bx,(0,\by);s)$,

\noi $L_{n-1,n}^0(s):\E_n\to\E_{n-1}^+$ with kernel $\ph((0,\bx),\by;s)$.

\noi(a) The operator $M_{n}(s)\,L_{n}(s)\,A_n:\E_{n-1}\to\E_n$ in Corollary 1.1 is equal to
\[L_{n,n-1}^0(s)\,(I-(\a-\b)\,L_{n-1}^0(s))\inv\,R_{n-1},
\footnote{The $I$ here is the identity operator on $\E_{n-1}^+$ while the $I$ in (\ref{Mn}) is the identity operator on $\E_n$.}\]
where $R_{n-1}:\E_{n-1}\to\E_{n-1}^+$ is the restriction operator.

\noi(b) The operator $M_{k}(s)\,L_{k}(s):\E_k\to\E_k$ in Corollaries 1.1 and 1.2 is equal to 
\[L_k(s)+(\a-\b)\,L_{k,k-1}^0(s)\,(I-(\a-\b)\,L_{k-1}^0(s))\inv\,
L_{k-1,k}^0(s).\]  
 
\sp

\noi{\bf Remark 1.2}. In the special case of TASEP when $p=1$ we have the simplification $(I-\a\,L_{n}^0(s))\inv=I+\a\,L_{n}^0(s-\a)$.  
\sp
 
In the case of SSEP ($p=q$), even for infinitly many particles, there are formulas for correlations that are no more complicated when both $\a$ and $\b$ are nonzero. For \linebreak$\bx=\{x_1,\ld,x_n\}\in\X_n$ we define $\Ps_n(\bx;t)$ to be the probability that sites $x_1,\ld,x_n$ are occupied at time $t$. We denote its Laplace transform by $\Psh_n(s)=\Psh_n(\bx;s)$
and introduce the vector functions
\[\Psh(s)=(\Psh_n(s))_{n\ge0},\ \ \ \Ps(0)=(\Ps_n(0))_{n\ge0}.\footnote{We define $\Ps_0(t)=1$, and so $\Psh_0(s)=s\inv$.}\]

Let the operators $L_n(s)$ and $A_n$, and the matrices $\d,\,L(s),\,A$ with operator entries be the same as before.\footnote{Now we define $\E_n=L^\iy(\X_n)$, and observe that by (\ref{psum}) and the symmetry of the kernel we have
$\sum_{\by\in\X_n}\p(\bx,\by;t)=1$
for each $\bx\in\X_n$. It follows that $L_n(s)$ is a bounded operator on this $\E_n$ with norm at most $({\rm Re}\,s)\inv$.}
\sp

\noi{\bf Theorem 2}. We have,
\[\Psh(s)=\Big(I+L(s)\,((\a+\b)\,\d-\a\,A)\Big)\inv\,L(s)\,\Ps(0).\]

We now set
\[M_n(s)=(I+(\a+\b)\,L_n(s)\,\d)\inv.\]

\noi{\bf Corollary 2.1}. For $n>0$,
\[\Psh_n(s)=\a\,M_n(s)L_n(s)A_n\,\Psh_{n-1}(s)+M_n(s)L_n(s)\,\Ps_n(0).\]

\noi{\bf Corollary 2.2}. In the case of Bernoulli initial condition with density $\r$, 
\[\Psh_n(s)=\a\,M_n(s)L_n(s)A_n\,\Psh_{n-1}(s)+s\inv\,\r^n\,M_n(s)\,1,\]
where ``1'' is the constant function on $\X_n$.

\noi{\bf Corollary 2.3}. When initially no sites are occupied,  
\[\Psh_n(s)=\a\,M_n(s)L_n(s)A_n\,\Psh_{n-1}(s).\]

The analogue of Remark 1.1 holds here.

From the abstract formulas, Theorems 1 and 2 and their corollaries, we derive some concrete results. 

Suppose, in ASEP, that at time zero there is a single particle at $y$. From Corollaries~1.1 and 1.2 we show:

\noi{\boldmath$\a=0$}: When $p>q$, with probability
\[1-{\b\,(q\inv-1)^{-y}\ov p-q+\b}\]
the particle is never ejected. When $p\le q$, with probability one the particle will eventually be ejected. The expected value of the time this occurs is infinite when $p=q$, and when $p<q$ it is 
\[{y+q/\b\ov q-p}.\]
\sp

\noi{\boldmath$\b=0$}: With probability one a second particle will eventually be injected. The expected value of the time at which this occurs is
\[{1\ov\a}+{\x_+(\a)^{-y}\ov q\,(\x_+(\a)-1)-\a},\]
where
\[\x_+(\a)={1\ov 2q}\,(\a+1+\sqrt{(\a+1)^2-4pq}\,)\,.\]
\sp

Combining Remark 1.2 with the determinant formula (\ref{TASEP}) for $\p(\bx,\by;t)$ in TASEP it is practical to compute some exact results. Denote by $\P_n(t)$  the probability that starting with no particles at time 0 there are exactly $n$ particles at time $t$. Then
\[\P_0(t)=e^{-\a t},\]
\[\P_1(t)={\a\ov1-\a}\left(t-{\a\ov1-\a}\right)e^{-\a t} +{\a^2\ov(1-\a)^2}e^{-t},\] 
\[\P_2(t) =\left({\a^2\ov2(1-\a)^2} t^2-{\a^2\ov(1-\a)^3}t+{\a^2\ov(1-\a)^2}\right)e^{-\a t}\]
\[+\left({\a^2\ov2(1-\a)^2}t^2-{\a(1-3\a+\a^2)\ov(1-\a)^3}+{1-2\a\ov(1-\a)^2}\right)e^{-t}  - e^{-(1+\a)t}.\]
These are for $\a\ne1$. When $\a=1$ there is no singularity; we take a limit and the formulas simplify.
\begin{figure}[h]
\bc \includegraphics{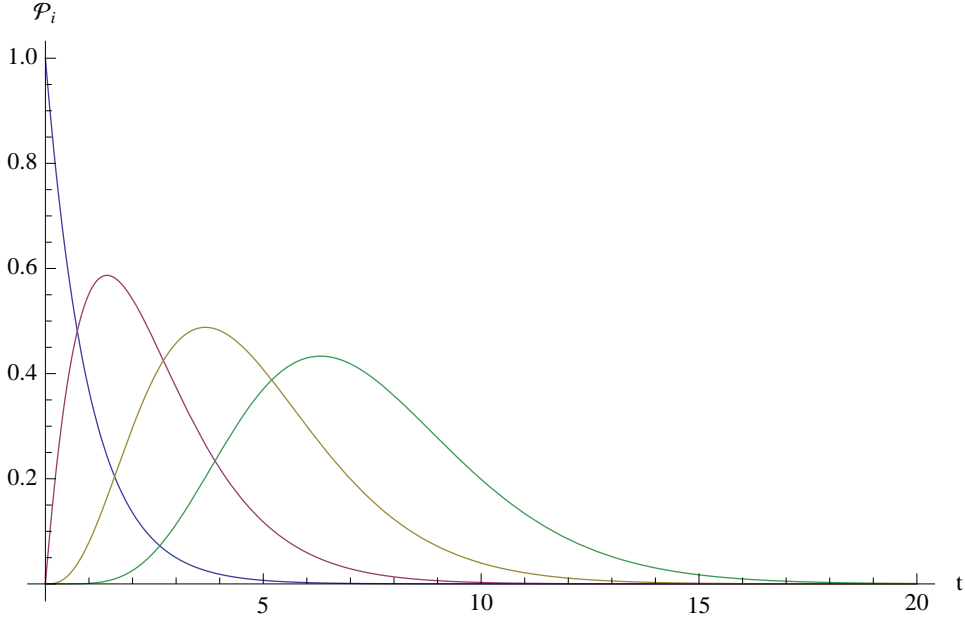}\ec
\caption{Plotted are the TASEP probabilities $\mathcal{P}_n(t)$  for $n=0,1,2, 3$ with $\alpha=1$.  Increasing $n$ moves the maximum to the right.}
\end{figure}

For SSEP we obtain the following consequence of Corollary 2.2. Define
\[\D N(t)=\sum_{x\ge0}\left( \eta_x(t)-\eta_x(0)\right),\]
the net number of particles that have entered the system at time $t$. (Which may be negative.) With Bernoulli initial condition when $\lb\e_x(0)\rb=\r$, we show for the expected value that
\[ \lb \D N(t)\rb \sim\sqrt{2\ov\pi}\,{\a-(\a+\b)\rho\ov\a+\b}\,t^{1/2}\ \ {\rm as}\  t\to\iy.\footnote{This result was obtained in \cite{K} in the case $\b=0$  and $\r=0$.}\]
By a laborious computation (not included) we can show that the second moment of $\D N(t)$ is finite. 

In the special case where initially there are no occupied sites, we use Corollary~2.3 and present a not completely rigorous (to say the least) argument that 
\[\lb \D N(t)^2\rb\sim{2\ov\pi}{\a^2\ov(\a+\b)^2}\,t\ \ \ {\rm as}\ t\to\iy.\]
Combining this with the first moment asymptotics when $\r=0$
we conclude that the variance of $\D N(t)/t^{1/2}$ tends to zero as $t\to\iy$.\footnote{As predicted in \cite{K} in the case $\b=0$.} The derivation is quite long, but the result with the precise constants came out so nicely in the end that we could not resist including it.

\bc{\bf III. Proofs of the results}\ec

\bc{\bf A. ASEP}\ec

\noi{\bf Proof of Theorem 1}: The probability $\p(\bx,\by;\,t)$ for $n$-particle ASEP on $\Z^+$ is the solution of the differential equation 
\[{d\ov dt}\,u(X;t)=\sum_{i=1}^N
\Big[ p\,u(x_i-1)\,(1-\d(x_i-x_{i-1}-1))
+q\,u(x_i+1)\,(1-\d(x_{i+1}-x_i-1))\]
\[ -p\,u(x_i)\,(1-\d(x_{i+1}-x_i-1))-q\,u(x_i)\,(1-\d(x_i-x_{i-1}-1))\Big]+[q\,u(x_1)-p\,u(x_1-1)]\,\d(x_1)\]
that satisfies the initial condition
\[u(\bx;\,0)=\d_\by(\bx).\]
(In the $i$th summand entry $i$ is displayed and entry $j$ is $x_j$ when $j\ne i$.)

If we denote by $\Q_n$ the operator given by the right side of the equation, then $\p(\bx,\by;\,t)$ is the kernel of $e^{t\Q_n}$. Thus the  equation is
\[{d\p_n\ov dt}=\Q_n\,\p_n(\bx;\,t),\]
where we have written $\p_n(\bx;\,t)$ for $\p(\bx,\by;\,t)$.

For open ASEP we write the corresponding probability as $P_n(\bx;\,t)$ (and do not specify any initial condition). The equation for $P_n$ is 
\[{d P_n\ov dt}(x_1,\ld,x_n:\,t)=\Q_n\,P_n(x_1,\ld,x_n:\,t)\]
\[+\a\,\d(x_1)\,P_{n-1}(x_2,\ld,x_n:\,t)
-\a\,(1-\d(x_1))\,P_n(x_1,\ld,x_n:\,t)\]
\[-\b\,\d(x_1)\,P_n(x_1,\ld,x_n:\,t)+\b(1-\d(x_1))\,P_{n+1}(0,x_1,\ld,x_n:\,t).\]

Define $\Ph_n(\bx;\,t)=e^{\a t}\,P_n(\bx;\,t)$. In terms of the operators $A_n$ and $B_n$ defined in (\ref{A}) and (\ref{B}), and $\d=\d(x_1)$, the equation for $\Ph_n$ becomes
\be{d\Ph_n(t) \ov dt}=\Q_n\,\Ph_n+(\a-\b)\,\d\,\Ph_n+\a\,A_n\,\Ph_{n-1}+\b\,B_n\,\Ph_{n+1}.\label{Pheq}\ee
The equation and initial condition are satisfied if
\[\Ph_n(t)=\int_0^t e^{(t-u)\Q_n}\,\Big((\a-\b)\,\d\,\Ph_n(u)+\a\,A_n\,\Ph_{n-1}(u)+\b\,B_n\,\Ph_{n+1}(u)\Big)\,du
+e^{t\Q_n}\,\Ph_n(0).\]

We use the fact that the Laplace transform of a convolution is the product of the Laplace transforms. Recall that the kernel of $L_n(s)$ is $\ph(\bx,\by;s)$, which is the Laplace transform of the kernel of $e^{t\Q_n}$. In other words, the Laplace transform of the operator $e^{t\Q_n}$ is $L_n(s)$.
So taking Laplace transforms in the last equation gives
\[\Phh_n(s)=L_n(s)\,\Big((\a-\b)\,\d\,\Phh_n(s)+\a\,A_n\,\Phh_{n-1}(s)+\b\,B_n\,\Phh_{n+1}(s)\Big)+L_n(s)\,\Ph_n(0).\]

If we now introduce the vector functions 
\[\Phh(s)=(\Phh_n(s))_{n\ge0},\ \ \ \Ph(0)=(\Ph_n(0))_{n\ge0},\]
and the operator matrices $\d,\,L(s),\,A,\, B$ defined earlier we see that the system may be written as
\[\Phh(s)=L(s)\,\Big((\a-\b)\,\d+\a\,A+\b\,B\Big)\,\Phh(s)+L(s)\,\Ph(0).\]
Since $\Phh(s)=\Phat(s-\a)$, this gives the statement of Theorem 1.
\sp

\noi{\bf Corollaries 1.1 and 1.2}: We write the operator inverse in (\ref{Th1}) as
\[\Big(I-(M(s)\,L(s)\,(\a\,A+\b\,B))\Big)\inv\,M(s)\,L(s),\]
where $M(s)$ is the diagonal matrix with operator entries $M_n(s)$. When $\b=0$ this equals
\be (I-\a\,M(s)\,L(s)\,A)\inv\,M(s)\,L(s).\label{b=0}\ee
The operator matrix $M(s)\,L(s)\,A$ consists of one subdiagonal, with $m,m-1$-entry $M_m(s)\,L_m(s)\,A_{m}$. Therefore the $n,n'$-entry of the inverse ($n'\le n$) is
\[\a^{n-n'}\,M_n(s)L_n(s)A_n\cd M_{n'+1}(s)L_{n'+1}(s)A_{n'+1},\]
where for $n=n'$ this equals $I$. So the $n,n'$-entry of (\ref{b=0}) is
\[\a^{n-n'}\,M_n(s)L_n(s)A_n\cd M_{n'+1}(s)L_{n'+1}(s)A_{n'+1}M_{n'}(s)L_{n'}(s).\]
Thus,
\be\Phat_n(s-\a)=\sum_{n'\le n}\a^{n-n'}\,M_n(s)L_n(s)A_n\cd 
M_{n'+1}(s)L_{n'+1}(s)A_{n'+1}\,M_{n'}(s)\,L_{n'}(s)\,P_{n'}(0),\label{Phatsum}\ee

When $\a=0$ the matrix $M(s)\,L(s)\,B$ consists of one superdiagonal, and we obtain similarly
\[\Phat_n(s)=\sum_{n'\ge n}\b^{n'-n}\,M_n(s)L_n(s)B_n\cd
M_{n'-1}(s)L_{n'-1}(s)B_{n'-1}\,M_{n'}(s)\,L_{n'}(s)\,P_{n'}(0).\]

If initially there are $k$ particles at $\by\in\X_k$ then in both cases $P_{n'}(0)$ is nonzero only for $n'=k$. The formulas become
\[\Phat_n(s-\a)=\a^{n-k}\,M_n(s)L_n(s)A_n\cd 
M_{k+1}(s)L_{k+1}(s)A_{k+1}\,M_{k}(s)\,L_{k}(s)\,\d_y,\]
\[\Phat_n(s)=\b^{k-n}\,M_n(s)L_n(s)B_n\cd
M_{k-1}(s)L_{k-1}(s)B_{k-1}\,M_{k}(s)\,L_{k}(s)\,\d_y,\]
and the corollaries follow.
\sp

\noi{\bf Remark 1.1}: For (a) we use the fact that
because $A_n$ has the factor~$\d$,   
\[M_n(s)L_n(s)A_n=(I-\a\,L_n(s)\,\d)\inv L_n(s)\,\d\,A_n
=L_n(s)\,\d\,(I-\a\,\d\,L_n(s)\,\d)\inv\,A_n.\]
Since $(A_n f)(x_1,\bx)=\d(x_1)\,(R_{n-1}f)(\bx)$ and $\bx\in\E_{n-1}^+$, we obtain statement (a) in different notation.
For (b) we use
\[M_k(s)=(I-(\a-\b)\,L_k(s)\d)\inv
=I+(\a-\b)\,L_k(s)\d\,(I-(\a-\b)\,L_k(s)\d)\inv\]
\be=I+(\a-\b)\,L_k(s)\d\,(I-(\a-\b)\,\d L_k(s)\d)\inv.\label{Mk}\ee
Thus,
\[M_k(s)L_k(s)=L_k(s)+(\a-\b)\,L_k(s)\d\,(I-(\a-\b)\,\d L_k(s)\d)\inv \d L_k(s),\]
and statement (b) follows. 
\sp

\noi{\bf Initially a single particle}:
Let $\P_1(y;t)=\sum_{x\ge0}P_1(x,y;t)$ denote the probability that, starting with one particle at $y$, we still have one particle at time~$t$. Denote its Laplace transform by $\LP_1(y;s)$.

We begin with the case $\b=0$, so $\P_1(y;t)$ is the probability that no new particle has been injected by time $t$. 

By Corollary 1.1 the Laplace transform $\widehat P_1(x,y;s-\a)$ is equal to $(M_1(s)L_1(s)\d_y)(x)$. By Remark 1.1(b) this equals
\[\ph(x,y;s)+\a\,\ph(x,0;s)\,(1-\a\,\ph(0,0;s))\inv\,\ph(0,y;s).\]
Then using $\sum_{x\ge0}\ph(x,y;s)=s\inv$ we obtain
\[\LP_1(y;s-\a)={1\ov s}\Big[1+{\a\,\ph(0,y;s)\ov 1-\a\,\ph(0,0;s)}\Big].\]

From formula (\ref{ASEP}) for $\p(x,y;t)$ in the case $n=1$ we compute that
\be\ph(0,y;s)={1\ov q}{\x_+(s)^{-y}\ov\x_+(s)-1},
\label{p0y}\ee
where 
\[\x_+(s)={1\ov 2q}\,(s+1+\sqrt{(s+1)^2-4pq});\]
this is the solution of $\ep(\x)=s$ with positive square root when  $s>0$. 
Thus
\be\LP_1(y;s-\a)={1\ov s}\left[1+{\a\,\x_+(s)^{-y}\ov q(\x_+(s)-1)-\a}\right].\label{N}\ee
The denominator in the brackets is nonzero for $s=\a$ (and positive for $s>\a$), from which we conclude that
\[\int_0^\iy\P_1(y;t)\,dt=\LP_1(y;0)={1\ov\a}+{\x_+(\a)^{-y}\ov q(\x_+(\a)-1)-\a}.\]
It follows that with probability one a second particle will eventually be injected, since $\P(y;t)\to0$ as $t\to\iy$, and the integral is the expected time at which it occurs.\footnote{If $T(y)$ denotes the time when a second particle is injected, then ${\rm Prob}\,(T(y)>t)=\P_1(y;t)$, from which the statement follows.}

This was for $\b=0$. For $\a=0$, $\P_1(y;t)$ is the probability that the particle has not been ejected by time $t$. We use Corollary 1.2 (and Remark 1.1), and formula (\ref{N}) is replaced by
\[\LP_1(y;s)={1\ov s}\left[1-{\b\,\x_+(s)^{-y}\ov q(\x_+(s)-1)+\b}\right].\]

We compute that as $s\to0$,
\[\LP_1(y;s)\sim {1\ov s}\left[1-{\b\,(q\inv-1)^{-y}\ov p-q+\b}\right]\ \ \ \textrm{if}\ p>q,\]
\[\LP_1(y;s)\sim {1\ov\sqrt{2s}}\,(2y+1/\b)\ \ {\rm if}\ p=q,\]
\[\LP_1(y;s)\to {y+q/\b\ov q-p}\ \ {\rm if}\ p<q.\]

Applying the Tauberian theorem we deduce
\[\lim_{t\to\iy}\P_1(y;t)=1-{\b\,(q\inv-1)^{-y}\ov p-q+\b}\ \ \ \textrm{if}\ p>q,\]
\[\P_1(y;t)\sim {1\ov\sqrt{2\pi}}\,\Big[2y+{1\ov\b}\Big]\,t^{-1/2}\ \ {\rm as}\ t\to\iy\ {\rm if}\ p=q,\]
\[\int_0^\iy \P_1(y;t)\,dt={y+q/\b\ov q-p}\ \ {\rm if}\ p<q.\]

When $p>q$ the limit on the first line is the probability that the particle is never ejected. If $p\le q$, then with probability 1 the particle will eventually be ejected. The expected value of the time at which it is ejected is infinite when $p=q$, by the second line, and is given by the integral on the last line when $p<q$.

\bc{\bf B. TASEP}\ec
  
\noi{\bf Remark 1.2}: To compute $(I-\a\,L_n^0(s))\inv$ we use the fact that for $k\ge1$ the kernel of $L_n^0(s)^k$ is
\[\sum_{\bz_1,\ld, \bz_{k-1}\in\X_n^+}\ph((0,\bx),(0,\bz_{1});s)\;
\ph((0,\bz_{1}),(0,\bz_{2});s)\cd\ph((0,\bz_{k-1}),(0,\by);s).\]
The summand is the Laplace transform of the $(k-1)$-fold convolution
\[\int_{u_1+\cd+u_{k-1}=t}\p((0,\bx),(0,\bz_{1});u_{1})\;
\p((0,\bz_{1}),(0,\bz_{2});u_{2})\cd \p((0,\bz_{n-1}),(0,\by);u_{k-1})\,d{\bf u}.\]
When $p=1$, if the left-most particle begins at 0 and ends at 0 then it was always at 0. It follows that  these probabilities have the semigroup property. Thus after summing over $\bz_1,\ld,\bz_{k-1}$ the integral becomes
\[\int_{u_1+\cd+u_{k-1}=t}\p((0,\bx),(0,\by),t)\,d{\bf u}={t^{k-1}\ov(k-1)!}\,\p((0,\bx),(0,\by),t).\]
Therefore the kernel of of $L_n^0(s)^k$ is
\[\int_0^\iy e^{-st}\,{t^{k-1}\ov(k-1)!}\,\p((0,\bx),(0,\by),t)\,dt.\]
It follows that the kernel of $\sum_{k=1}^\iy (\a\,L_n^0(s))^k$ is
\[\a\,\int_0^\iy e^{-st}\,e^{\a t}\,\p((0,\bx),(0,\by),t)\,dt
=\a\,\ph((0,\bx),(0,\by);s-\a).\]
This gives $(I-\a\,L_n^0(s))\inv=1+\a\,L_n^0(s-\a)$.

Once we have this result (and Corollary 1.2, Remark 1.1, and the determinant formula (\ref{TASEP}) for $\p(\bx,\by;t)$), we have all the ingredients necessary to compute $\P_n(t)$ for small $n$.  This is the probability that starting with no particles at time 0 there are exactly $n$ particles at time $t$, which equals $\sum_{\bx\in\X_n}\P_n(\bx;t)$. The computation of products of kernels involves summing geometric series.

\bc{\bf C. SSEP}\ec

\noi{\bf Proof of Theorem 2}: A state of the system is a  function $\eta:\Z^+\to \{0,1\}$ where $\eta_x=1$ means site~$x$ is occupied and $\eta_x=0$ means site~$x$ is not occupied. Recall that we defined $\Ps_n(x_1,x_2,\ld,x_n;\,t)$ as the probability that sites $x_1,\ld,x_n$ are occupied at time $t$. Thus,
\[\Ps_n(x_1,x_2,\ld,x_n;\,t)=\lb\e_{x_1}(t)\cd\e_{x_n}(t)\rb.\]

The Markov generator $\L$ of ASEP on $\Z^+$ with an open boundary at zero is given \cite{SW} by 
\[\L f(\e)=\a\,(1-\e_0)\left(f(\e^0)-f(\e)\right) +\b\,\e_0 \left(f(\e^0)-f(\e)\right)\]
\be+\sum_{k=0}^\iy \left[p\,\e_k\,(1-\e_{k+1})+ q\,(1-\e_k)\,\e_{k+1}\right] \left(f(\e^{k,k+1})-f(\e)\right).\label{generator}\ee
Here $f$ is an $\R$-valued function that depends on only finitely many sites, and 
\[(\e^k)_x=\left\{\begin{array}{ll} 1-\e_x & \textrm{if}\>\> x=k \\
\e_x& \textrm{if}\>\> x\neq k,\end{array}\right.\ \ \ \ \left(\e^{k,k+1}\right)_x=\left\{\begin{array}{ll} \e_{k+1} & \textrm{if}\>\> x=k \\
  \e_k & \textrm{if}\>\> x=k+1 \\
  \e_x & \textrm{if}\>\> x\neq k, k+1.\end{array}\right.\]

For SSEP the sum in (\ref{generator}) equals $1/2$ times
\[\sum_{k=0}^\iy [\e_k\,(1-\e_{k+1})+(1-\e_k)\,\e_{k+1}]\,(f(\e^{k,k+1})-f(\e)).\]
When $\e_k=\e_{k+1}$ the first factor equals 0; otherwise it equals 1. Since the second factor is zero when $\e_k=\e_{k+1}$, we can ignore the first factor, and we get
\[\sum_{k=0}^\iy (f(\e^{k,k+1})-f(\e)).\]

For the correlations we are interested in $f(\e)=\e_{x_1}\cd\e_{x_n}$, so the $k$th summand equals zero unless either $k=x_i$ for some $i$ or $k=x_i-1$ for some $i$. (We'll see that these cannot both happen for a nonzero summand.)

Suppose first that $k=x_i$.
If $x_{i+1}=x_i+1$ then $k+1=x_{i+1}$ and the substitution $\e\to\e^{k,k+1}$ applied to $f(\e)$ just interchanges $\e_{x_i}$ and $\e_{x_{i+1}}$. Therefore the $k$th summand is zero. It follows that for a nonzero summand we must have $x_{i+1}>x_i+1=k+1$, and the substitution $\e\to\e^{k,k+1}$ only affects the $i$th factor in $f(\e)$. Therefore the summand equals
\be\Big(\e_{x_1}\cd\e_{x_{i-1}}\,\e_{x_i+1}\,\e_{x_{i+1}}\cd\e_{x_n}-
\e_{x_1}\cd\e_{x_{i-1}}\,\e_{x_i}\,\e_{x_{i+1}}\cd\e_{x_n}\Big)\,(1-\d(x_{i+1}-x_i-1)).\label{first}\ee
Similarly if $k=x_i-1$ the summand equals
\be\Big(\e_{x_1}\cd\e_{x_{i-1}}\,\e_{x_i-1}\,\e_{x_{i+1}}\cd\e_{x_n}-
\e_{x_1}\cd\e_{x_{i-1}}\,\e_{x_i}\,\e_{x_{i+1}}\cd\e_{x_n}\Big)\,(1-\d(x_i-x_{i-1}-1)).\label{second}\ee
(This is to be multiplied by $1-\d(x_1)$ when $i=1$.)

If $k=x_i$ for (\ref{first}) and $k=x_{i'}-1$ for (\ref{second}) then $i'=i+1$, so $x_{i+1}=x_{i'}=x_i+1$ and (\ref{first}) zero, and $x_{i'}=x_i+1=x_{i'-1}+1$ so (\ref{second}) is zero. Thus for the $k$th summand to be nonzero either $k=x_i$ for some $i$ or $k=x_i-1$ for some $i$, but not both.

It follows that for SSEP the expected value of the sum in 
(\ref{generator}) is equal to $1/2$ times the sum over $i$ of the expected values of the sum of (\ref{first}) and (\ref{second}). This equals $\Q_n(\lb\e_{x_1}\cd\e_{x_n}\rb)$.

Adding what we get from the first two terms of (\ref{generator}) we find that the differential equation for $\Ps_n=\Ps_n(t)$ is
\[{d\ov dt}\Ps_n(x_1,x_2,\ld,x_n)=\Q_n\,\Ps_n(x_1,x_2,\ld,x_n)\]
\[-(\a+\b)\,\Ps_{n}(x_1,,x_2,\ld,x_n)\,\d(x_1)+\a\,\Ps_{n-1}(x_2,\ld,x_n)\,\d(x_1).\footnote{Observe that the constant functions $\Ps_n=(\a/(\a+\b))^n$ satisfy the equations. Thus the Bernoulli measure with density $\r=\a/(\a+\b)$ is stationary. This also follows from the results of \cite{SW} in which stationary measures were determined for general ASEP on $\Z^+$.}\]
Except for the change $\a-\b\to-(\a+\b)$, this is (\ref{Pheq}) without the $B_n$ terms. Therefore to complete the proof of Theorem 2 we need only make this change as we go through the rest of the proof of Theorem 1, which we need not do.
\sp
\noi{\bf Corollaries 2.1 and 2.2}: Just as (\ref{Phatsum}) is obtained from Theorem 1, we obtain now for $n>0$,
\[\Psh_n(s)=\sum_{k=0}^n\a^{n-k}\,M_n(s)L_n(s)A_n\cd M_{k+1}(s)L_{k+1}(s)A_{k+1}\,M_k(s)L_k(s)\,\Ps_k(0).\]
Corollary 2.1 follows. For Corollary 2.2 we have $\Ps_n(0)=\r^n$, and we use the fact that $L_n(1)=\sum_{\by\in\X_n}\ph(\bx,\by;\,s)=s\inv$. Corollary 2.3 is the case $\r=0$ of Corollary 2.2.
\sp

\noi{\bf \boldmath$\D N(t)$, the net number of particles that entered the system}: We assume that we have SSEP with Bernoulli initial condition. Corollary 2.2 when $n=1$ gives, for the Laplace transform of $\lb\e_x(t)\rb$,
\[\wh{\lb\e_x\rb}(s)=\Psh_1(x;s)=s\inv\,(\r\,M_1(s)\,1+\a\,M_1(s)L_1(s)A_1\,1).\]
By Remark 1.1(a) applied here, we have\footnote{Since the definition of $M_k(s)$ is different now we must replace $\a-\b$ by $-(\a+\b)$ when using the remark.}
\be M_1(s)L_1(s)A_1\,1={\ph(x,0;s)\ov 1+\g\,\ph(0,0;s)},\label{ev1}\ee
where we set
\[\g=\a+\b.\]
Similarly, from (\ref{Mk}) we obtain
\[ M_1(s)\,1=1-\g{\ph(x,0;s)\ov 1+\g\,\ph(0,0;s)}.\]
Combining the two gives  
\[\wh{\lb\e_x\rb}(s)={1\ov s}\Big[\r+{\a-\g\,\r\ov 1+\g\,\ph(0,0;s)}\,
\ph(x,0;s)\Big].\]
Subtracting $\r/s$ from both sides and summing over $x\ge0$ we get for the Laplace transform of $\lb \D N(t)\rb$ 
\[\wh{\lb \D N(t)\rb}(s)={1\ov s^2}\,{\a-\g\r
\ov 1+\g\,\ph(0,0;s)}.\]

{}From (\ref{p0y}) we have for SSEP
\[\ph(0,0;s)={2\ov s+\sqrt{s^2+2s}}\sim\sqrt{{2\ov s}}\ \ {\rm as}\ s\to0.\]
Hence
\[\wh{\lb \D N(t)\rb}(s)\sim{\a-\g\,\r\ov\sqrt2\,\g}\,s^{-3/2}\ \ {\rm as}\ s\to0.\]
By the Tauberian theeorem this implies
\[\lb \D N(t)\rb\sim\sqrt{{2\ov\pi}}\,{\a-\g\,\r\ov \g}\,t^{1/2}\ \ {\rm as}\ t\to\iy.\]\sp

As stated in the last section the second moment $\lb \D N(t)^2\rb$ is finite. To show this we use that the second moment is equal to
\[\lim_{N\to\iy}\sum_{x_1,x_2<N}\lb(\e_{x_1}(t)-\e_{x_1}(0))\,(\e_{x_2}(t)-\e_{x_2}(0))\rb.\]
We can show that the sum is a polynomial of degree two in $\r$, and that each of the three coefficients of the powers of $\r$ has a limit as $N\to\iy$. The argument is quite involved, and we do not include it. 
\sp

\noi{\bf The second moment of \boldmath$\D N(t)$ when \boldmath$\r=0$}: Now it is mainly a question of determining the asymptotics of $\sum_{x_1<x_2}\lb\e_{x_1}\e_{x_2}\rb$. The Laplace transform of $\Ps_2(\bx;t)=\lb\e_{x_1}\,\e_{x_2}\rb$
is given by Corollary 2.2 as
\[\Psh_2(\bx;s)=\a^2\,s^{-1}\,M_2(s)L_2(s)A_2M_1(s)L_1(s)A_1\,1.\]

We computed $M_1(s)L_1(s)A_1\,1$ in (\ref{ev1}). 
Combining this with Remark 1.1(a) gives
\[\Psh_2(\bx;s)=\a^2\,s\inv\,L_{2,1}^0(s)\,(I+\g\,L_1^0(s))\inv
{\ph_0\ov 1+\g\,\ph_0(0)},\] 
where $\ph_0$ is the function $x\to\ph(x,0;s)$.
   
To obtain the Laplace transform $\sum_{x_1<x_2}\wh{\lb\e_{x_1}\e_{x_2}\rb}$ 
we sum over $\bx\in\X_2$. If we recall that the kernel of $L_{2,1}^0(s)$ is $\ph(\bx,(0,y))$ and that the sum of this over $\bx\in\X_2$ is $s\inv$ we see that the desired sum is the inner product of the remaining function with the constant function $s\inv$. Thus, 
\be\sum_{x_1<x_2}\wh{\lb\e_{x_1}\e_{x_2}\rb}={\a^2s^{-2}\ov 1+\g\,\ph_0(0)}
((I+\g\,L_1^0(s))\inv\ph_0,1).\label{2nd}\ee
(Here we use the fact that $\ph_0$ belongs to $L^1$ and that $L_1^0(s)$ is a bounded operator on this space.)

What follows is not rigorous. We want to rescale $(I+\g\,L_1^0(s))\inv$ as $s\to0$, and we refer to formula (\ref{ASEP}) given in the Appendix for $\p(\bx,\by;t)$. After taking Laplace transforms in the case $n=2$ we find that 
$L_1^0(s)$ has kernel
\be L_1^0(x,y;s)
=\sum_{\s\in\B_2} {1\ov(2\pi i)^2}\int_{\C_R}\int_{\C_R} A_\s{\xs{2}^x\,\x_2^{-y}\ov s-\ep(\x_1)-\ep(\x_2)}{d\x_1\,d\x_2\ov\x_1\,\x_2},\label{kernel}\ee
where $\C_R$ is a circle with radius $R$ with $R$ large. (Some integrals are taken over two pairs of different contours and the results averaged.) To begin with, $s$ is so large that taking Laplace transforms under the integral sign is valid. 

If we ignore the poles of the $A_\s$ we can move both contours to the unit circle $\C_1$. Then the range of $\ep(\x_1)+\ep(\x_2)$ is $[-4,0]$, so we may take any $s>0$. As $s\to0$ the main contribution comes from a neighborhood of $\x_1=\x_2=1$ because the denominator vanishes there when $s=0$. If we set $\x_1=e^{iv_1},\,\x_2=e^{iv_2}$ then the integral with its factor becomes to first order
\[{1\ov4\pi^2}\int_\R\int_\R A_\s(e^{iv_1},e^{iv_2})\,{e^{i\,(\vs{2}\,x-v_2\, y)}\ov s+(v_1^2+v_2^2)/2}\,dv_1\,dv_2\]
\[={1\ov2\pi^2}\int_\R\int_\R A_\s(e^{i\sq\, v_1},e^{i\sq\, v_2})\,{e^{i\,\sq\,(\vs{2}\,x-v_2\, y)}\ov 1+v_1^2+v_2^2}\,dv_1\,dv_2.\]

This becomes, after the scaling $x\to x/\sq,\,y\to y/\sq$,
\[{1\ov2\pi^2\sq}\int_\R\int_\R A_\s(e^{i\sq\, v_1},e^{i\sq\, v_2})\,{e^{i\,(\vs{2}\,x-v_2\, y)}\ov 1+v_1^2+v_2^2}\,dv_1\,dv_2,\]
which acts on functions on $\R^+$.

Each $A_\s$ has absolute value 1 on $\C_1\times\C_1$, and each $A_\s(e^{i\sq v_1},e^{i\sq v_2})$ has limit~1 as $s\to0$ except when $v_1=v_2$ (mod $2\pi$). Thus we replace the above by the approximation
\[{1\ov2\pi^2\sq}\int_\R\int_\R {e^{i\,(\vs{2}\,x-v_2\, y)}\ov 1+v_1^2+v_2^2}\,dv_1\,dv_2.\]
This depends only on $\s(2)$. If we use the fact that the denominator is even in each $v_i$, and that $v_{-i}=-v_i$, we see that the sum over $\s\in\B_2$ of the integrals equals
\[{1\ov \pi^2\sq}\int_\R\int_\R {e^{i\,v_2\,(x-y)}+e^{i\,v_2\,(x+y)}\ov 1+v_1^2+v_2^2}\,dv_1\,dv_2+
{2\ov \pi^2\sq}\int_\R\int_\R {e^{i\,(v_1\,x-v_2\,y)}\ov 1+v_1^2+v_2^2}\,dv_1\,dv_2.\]
In the first integral we integrate first with respect to $v_1$, and we obtain
\be{1\ov \pi\sq}\int_\R {e^{i\,v\,(x-y)}+e^{i\,v\,(x+y)}\ov \sqrt{1+v^2}}\,dv+{2\ov \pi^2\sq}\int_\R\int_\R 
{e^{i\,(v_1\,x-v_2\,y)}\ov 1+v_1^2+v_2^2}\,dv_1\,dv_2.\footnote{The integrals are $K_0$ Bessel functions, but this fact is not useful.\label{factor}}\label{integrals}\ee

These are of the order $1/\sqrt{s}$ as $s\to0$. Now we indicate why the contributions from the poles of the $A_\s$, when we shrink the contours, are of lower order.

Consider the permutations $(\pm2\ 1)$, when $A_\s=S(\x_2,\x_1)$. (Times $\x_1\inv$ when $\s=(-2\ 1)$; this has no effect on what follows.) With the $\x_2$-integration over $\C_R$, we shrink the $\x_1$-contour to $\C_1$. Then when we shrink the $\x_2$-contour we pass the pole at $\x_2=2-\x_1\inv$ for all $\x_1\in\C_1$. The residue is a constant times 
\be\int_{\C_1}\left({\x-1\ov\x}\right)^2\,{\x^{x}\,(2-\x\inv)^{-y-1}
\ov s-{(\x-1)^2\ov2\x-1}}\,{d\x\ov\x},\label{res1}\ee
where we replaced $\x_1$ by $\x$. With either branch of $\x^{1/2}$ we may write
\[{(\x-1)^2\ov2\x-1}={(\x^{1/2}-\x^{-1/2})^2\ov2-\x\inv}.\]
On $\C_1$ the numerator is negative real (except when $\x=1$) while the denominator lies in the right half-plane. Thus the quotient lies in the left half-plane (except when $\x=1$) so we may take any $s>0$ in the integral. Since again the main contribution comes from near $\x=1$, we set $\x=e^{iv}$ and make the replacements $2-\x\inv\to2-(1-iv)=1+iv\to e^{iv}$, and we get
\[\int_\R{v^2\ov s+v^2}\,e^{iv(x-y)}\,dv=\sqrt{s}\int_\R {v^2\ov 1+v^2}\,e^{i\sqrt{s}\,v(x-y)}\,dv.\]
After the scaling this becomes independent of $s$. 

The factor $(\x-1)^2$ in the integrand in (\ref{res1}) was important. It also appears in the residues for the other integrals, which also become independent of $s$ by similar computations. We omit the details.
 
Thus when we scale the contributions from the poles of the $A_\s$ are independent of $s$, and so of lower order than the main terms (\ref{integrals}). 

Set
\[ J_1(x,y)={1\ov 2\pi}\int_\R {e^{i\,v\,(x-y)}+e^{i\,v\,(x+y)}\ov \sqrt{1+v^2}}\,dv,\]
\[ J_2(x,y)={1\ov\pi^2}\,\int_\R\int_\R 
{e^{i\,(v_1\,x-v_2\,y)}\ov 1+v_1^2+v_2^2}\,dv_1\,dv_2.\]
We showed that $1+\g L_1^0(s)$, when scaled, is equal to $2\g/\sq\,(J_1+J_2)$ plus an operator independent of $s$. Therefore we presume that to a first approximation the scaled operator $(1+\g L_1^0(s))\inv$ is equal to 
$(2\g)\inv\sq\,(J_1+J_2)\inv$. Also, we see from (\ref{p0y}) and the symmetry of $\ph(\bx,\by;s)$ that 
\[\ph_0(x)=\ph(x,0;s)\sim\sqrt{2\ov s}\,e^{-\sq\, x}\ \ \ 
{\rm as}\ s\to0.\]  
This gives the {\bf Conjecture}
\be((I+\g L_1^0(s))\inv\,\ph_0,1)\sim{1\ov \g\sq}((J_1+J_2)\inv\,e^{-x},1)\ \ \ {\rm as}\ s\to0.\label{conjecture}\ee

We shall show that the inner product equals $2/\pi$. Assume this, and the conjecture, for now. Then from (\ref{2nd}) and the asymptotics $\ph(0,0;s)\sim\sqrt{2/s}$ as $s\to0$ we find that  
\[\sum_{x_1<x_2}\wh{\lb\e_{x_1}\e_{x_2}\rb}\sim{1\ov\pi}{\a^2\ov \g^2}\,s^{-2}\ \ \ {\rm as}\ s\to0.\]
Therefore by the Tauberian theorem,
\[\sum_{x_1<x_2}\lb\e_{x_1}\,\e_{x_2}\rb\sim{1\ov\pi}{\a^2\ov\g^2}\,t\ \ \ {\rm as}\ t\to\iy.\]
Since
\[\lb \D N(t)^2\rb=\sum_{x_1,x_2\ge0}\lb\e_{x_1}\,\e_{x_2}\rb
=2\,\sum_{x_1<x_2}\lb\e_{x_1}\,\e_{x_2}\rb+\lb \D N(t)\rb,\]
and $\D N(t)=o(t)$, we get the result stated in Section II,
\[\lb \D N(t)^2\rb\sim{2\ov\pi}{\a^2\ov\g^2}\,t\ \ \ {\rm as}\ t\to\iy.\]

Now we show that the inner product in (\ref{conjecture}) equals $2/\pi$. We first use an observation about Wiener-Hopf plus Hankel operators, of which $J_1$ is one. (Recall that our operators act on functions on $\R^+$.) Suppose we have such an operator with kernel $J(x-y)+J(x+y)$, so that the result of its action on a function $f$ on $\R^+$ is
\[(Jf)(x)=\int_0^\iy(J(x-y)+J(x+y))\,f(y)\,dy.\]
If we extend $f$ to an even function on $\R$, then this equals
\[\int_{-\iy}^\iy J(x-y)\,f(y)\,dy.\]
So the operator becomes simply convolution by $J$.

Also, it follows from the fact that $J_2(x,y)$ is even in $y$ that if $f$, defined on~$\R^+$, is extended to be even on $\R$ then
\[(J_2\,f)(x)={1\ov2}\int_{-\iy}^\iy J_2(x,y)\,f(y)\,dy.\]

To recapitulate, define
\[K_1(x,y)={1\ov2\pi}\int_\R {e^{i\,v\,(x-y)}\ov \sqrt{1+v^2}}\,dv,\]
\[K_2(x,y)={1\ov2\pi^2}\int_\R\int_\R 
{e^{i\,(v_1\,x-v_2\,y)}\ov 1+v_1^2+v_2^2}\,dv_1\,dv_2,\]
acting on functions on $\R$. Then on $\R^+$ we have $(J_1+J_2)f=(K_1+K_2)f$, where the $f$ on the right is the even extension of the $f$ on the left. 

It follows from this that for an even function $f$ on $\R$, $(J_1+J_2)\inv f$ is the restriction to $\R^+$ of $(K_1+K_2)\inv f$.\footnote{This uses that $K_1+K_2$ commutes with the operator $f(x)\to f(-x)$.} Therefore the inner product in (\ref{conjecture}), which is over $\R^+$, is equal to 
\be{1\ov2}((K_1+K_2)\inv\,e^{-|x|},1),\label{ip}\ee
where this inner product is over $\R$.

Because of the forms of the kernels of $K_1$ and $K_2$, the operators simplify when we conjugate with the Fourier transform. The operator $K_1$ becomes $\Kh_1$, which is multiplication by the function $1/\sqrt{1+v^2}$, and the operator $K_2$ becomes $\Kh_2$, which has kernel 
\[\Kh_2(u,v)={1\ov\pi}{1\ov 1+u^2+v^2}.\]

Since the Fourier transform of $1$ with factor $1/2\pi$ outside the integral is 
$\d_0$ and the Fourier transform of $e^{-|x|}$ with factor 1 is 
$2/(1+v^2),$ we see that (\ref{ip}) (which equals the inner product in (\ref{conjecture})) equals 
\[((\Kh_1+\Kh_2)\inv(1+v^2)\inv,\d_0).\]
If we use 
\[(\Kh_1+\Kh_2)\inv=\Kh_1^{-1/2}\,
(1+\Kh_1^{-1/2}\,\Kh_2\,\Kh_1^{-1/2})\inv\,\Kh_1^{-1/2}\]
and the fact that $\Kh_1^{-1/2}\d_0=\d_0$, we see that the above equals
\be((I+\Kh)\inv\,\ps,\d_0),\label{ph}\ee
where  
\[\Kh(u,v)={1\ov\pi}{(1+u^2)^{1/4}\,(1+v^2)^{1/4}\ov 1+u^2+v^2},\ \ \ \ \ps(v)=(1+v^2)^{-3/4}.\]

If we conjugate with the unitary operator $f(u)\to (\cosh x)^{1/2}f(\sinh x)$
we find that (\ref{ph}) equals $((I+\Kt)\inv\,\pst,\d_0)$, where
\[\Kt(x,y)=(\cosh x)^{1/2}\,\Kh(\sinh x,\sinh y)\,(\cosh y)^{1/2}=
{1\ov\pi}{\cosh x\,\cosh y\ov1+\sinh^2x+\sinh^2y},\]
\[\pst(x)=(\cosh x)^{1/2}\,\ps(\sinh x)=\sech x.\]
{}From
\[\cosh x\,\cosh y={1\ov2}(\cosh(x+y)+\cosh(x-y)),\]
\[1+\sinh^2x+\sinh^2y=\cosh(x+y)\,\cosh(x-y)\]
one sees that
\[\Kt(x,y)={1\ov2\pi}(\sech(x-y)+\sech(x+y)).\]
Now $\pst$ is even, and when $\Kt$ is restricted to the space of even functions it equals the operator with convolution kernel $\sech(x-y)/\pi$. Conjugating with the Fourier transform, this operator becomes multiplication by
$\sech(\pi \x/2)$. The Fourier transform of $\sech x$ with factor $1/2\pi$ is $(1/2)\,\sech(\pi \x/2)$, and the Fourier transform of $\d_0$ with factor 1 is 1. Therefore
\[((I+\Kt)\inv\,\pst,\d_0)={1\ov2}\int_{-\iy}^\iy {1\ov 1+\sech(\pi \x/2)}\,\sech(\pi \x/2)\,d\x={2\ov\pi}.\]
 Thus the inner product in (\ref{conjecture}) is $2/\pi$ as claimed.

\bc{\bf IV. Appendix -- The formulas for \boldmath$\p(\bx,\by;t)$}\ec

The Weyl group $\B_n$ is the group of {\it signed permutations}, functions $\s:[1,\,n]\to[-n,\,-1]\cup[1,\,n]$ such that $|\s|\in\S_n$. An {\it inversion} in $\B_n$ is a pair $(\pm\s(i),\s(j))$ with $i<j$ such that $\pm\s(i)>\s(j)$.
We write $\t=p/q$. 

We define 
\[S(\x,\x')=-{p+q\,\x\x'-\x\ov p+q\,\x\x'-\x'},\ \ \ r(\x):={\x-1\ov 1-\t\,\x\inv},\ \ \ \ep(\x)=p\,\x\inv+q\,\x-1,\]
and then define
\[A_\s=\prod_{\s(i)<0}\,r(\xs{i})\,\times\,\prod\{S(\x_{a},\,\x_{b}):
(a,b)\ \textrm{is an inversion in}\ \B_n\},\]
with the convention $\x_{-a}=\t/\x_a$. 

The formula, valid when $q\ne0$, is
\be \p(\bx,\by;t)={1\ov n!}\sum_{\s\in\B_n} {1\ov(2\pi i)^n}\int\cd\int A_\s(\x)\,\prod_i\Big(\xs{i}^{x_i}\;\x_i^{-y_i-1}e^{\ep(\x_i)\,t}\Big)\,d\x_1\cd d\x_n,\label{ASEP}\ee
where $\bx=\{x_1,\ld,x_n\}$ and  $\by=\{y_1,\ld,x_n\}$.
The domain of integration is
\[\bigcup_{\m\in\S_n}\,\C_{\m(1)}\times\cd\times\C_{\m(n)},\]
where the $\C_a$ are circles with center $1/2q$ and distinct radii $R_a$. The $R_a$ should be so large that $S(\x,\x')$ is analytic for $\x,\,\x'$ on and outside $\C_a$.\footnote{We cannot simply take $\C\times\cd\times\C$ with $\C$ a circle with large radius, because then there would be nonintegrable singularities of $S(\x_a,\x_b)$ on the contour when $a>0,\,b<0$. However by taking the $R_a\to R$ we can interpret each integral as a symmetric distribution supported on $\C\times\cd\times\C$ applied to the product in the integrand.}

In the special case $n=1$,
\[\p(x,y;t)={1\ov 2\pi i}\int_{\C_R}\Big[\x^{x-y-1}+{\t-\x\ov1-\x}\,\t^x\,\x^{-x-y-2}\Big]\,e^{\ep(\x)\,t}\,d\x,\]
where $\C_R$ is a circle around 0 of radius $R>1$. 

The formulas do not hold for $p=1$ TASEP on $\Z^+$. But then the probability is the same as for TASEP on $\Z$, and we have then \cite{S} (or \cite[p. 820]{tw1})
\be \p(\bx,\by;t) =\det\left(\int_{\C_r} (1-\x)^{j-i} \x^{x_i-y_j-1} e^{t\ep(\x)}\, d\x\right), \label{TASEP}\ee
where $\C_r$ is a circle with center 0 and radius $r<1$.

\begin{center}{\bf Acknowledgments}\end{center}

We thank Lauren Williams for elaborating for us some of the results of \cite{SW}.  

This work was supported by the National Science Foundation through grants DMS-1207995 (first author) and DMS-0854934 (second author).

\end{document}